\documentclass[10pt]{article}
    \usepackage{amsmath,amsfonts,amsthm,amssymb,amscd,enumerate, url, anysize, hyperref}
   \usepackage[pdftex]{graphicx,color}     
\marginsize{1in}{1in}{1in}{1in}

\begin{document}

\title{What is... a $\mathrm{G}_2$-manifold?}

\author{Spiro Karigiannis \\
{\it Department of Pure Mathematics, University of Waterloo} \\
\tt{karigiannis@math.uwaterloo.ca} }

\maketitle

A $\mathrm{G}_2$-manifold is a Riemannian manifold whose holonomy group is contained in the exceptional Lie group $\mathrm{G}_2$. In addition to explaining this definition and describing some of the basic properties of $\mathrm{G}_2$-manifolds, we will discuss their similarities and differences to K\"ahler manifolds in general and to Calabi-Yau manifolds in particular.
 
The \emph{holonomy group} of a Riemannian manifold is a compact Lie group which in some sense provides a global measure of the local curvature of the manifold. If we assume certain nice conditions on the manifold and its metric, then, of the five exceptional Lie groups, only $\mathrm{G}_2$ can arise as such a holonomy group. Berger's classification in the 1950's could not rule them out, but it was generally believed that such metrics could not exist. Then in 1987 Robert Bryant successfully demonstrated the existence of \emph{local} examples. Two years later, Bryant and Simon Salamon found the first complete, non-compact examples of such metrics, on total spaces of certain vector bundles, using symmetry methods. Since then physicists have found many examples of non-compact holonomy $\mathrm{G}_2$ metrics with symmetry. Finally, in 1994 Dominic Joyce caused great surprise by proving the existence of hundreds of \emph{compact} examples. His proof is non-constructive, using hard analysis involving the existence and uniqueness of solutions of a non-linear elliptic equation, much as Yau's solution of the Calabi conjecture gives a non-constructive proof of the existence and uniqueness of Calabi-Yau metrics (holonomy $\mathrm{SU}(m)$ metrics) on K\"ahler manifolds satisfying certain conditions. In 2000 Alexei Kovalev found a different construction of compact manifolds with $\mathrm{G}_2$ holonomy that produced several hundred more non-explicit examples. These two are the only known compact constructions to date. An excellent survey of $\mathrm{G}_2$ geometry and some of the compact examples is~\cite{J}.

In terms of Riemannian holonomy, the aspect of the group $\mathrm{G}_2$ that is important is not that it is one of the five exceptional Lie groups, but rather that it is the automorphism group of the octonions $\mathbb O$, an $8$-dimensional \emph{non-associative} real division algebra. The octonions come equipped with a positive definite inner product, and the span of the identity element $1$ is called the \emph{real} octonions while its orthogonal complement is called the \emph{imaginary} octonions $\mathrm{Im} (\mathbb O) \cong \mathbb R^7$. This is entirely analogous to the quaternions $\mathbb H$, except that the non-associativity introduces some complications. This analogy allows us to define a \emph{cross product} on $\mathbb R^7$ as follows. Let $u, v \in \mathbb R^7 \cong \mathrm{Im} (\mathbb O)$ and define $u \times v = \mathrm{Im}(uv)$, where $uv$ denotes the octonion product. (In fact the real part of $uv$ is equal to $-\langle u, v \rangle$, just as it is for quaternions, where $\langle \cdot, \cdot \rangle$ denotes the Euclidean inner product.) This cross product satisfies the following relations:
\begin{equation*}
u \times v = - v \times u, \qquad \qquad \langle u \times v , u \rangle = 0, \qquad \qquad {|| u\times v||}^2 = {|| u \wedge v ||}^2, 
\end{equation*}
exactly like the cross product on $\mathbb R^3 \cong \mathrm{Im} (\mathbb H)$. However, there is a difference. Unlike the cross product in $\mathbb R^3$, the following expression is \emph{not} zero:
\begin{equation*}
u \times (v \times w) + \langle u, v \rangle w - \langle u, w \rangle v
\end{equation*}
but is instead a measure of the failure of the associativity: $(uv)w - u(vw) \neq 0$. We note that on $\mathbb R^7$ we can define a $3$-form (a totally skew-symmetric trilinear form) using the cross product as follows: $\varphi(u,v,w) = \langle u \times v, w \rangle$. For reasons we will not address here, this form is called the \emph{associative} $3$-form.

A $7$-dimensional smooth manifold $M$ is said to admit a $\mathrm{G}_2$-structure if there is a reduction of the structure group of its frame bundle from $\mathrm{GL}(7, \mathbb R)$ to the group $\mathrm{G}_2$, viewed naturally as a subgroup of $\mathrm{SO}(7)$. This implies that a $\mathrm{G}_2$-structure determines a Riemannian metric and an orientation. In fact, on a manifold with $\mathrm{G}_2$-structure, there exists a ``non-degenerate'' $3$-form $\varphi$ for which, at any point $p$ on $M$, there exist local coordinates near $p$ such that, at $p$, the form $\varphi$ is exactly the associative $3$-form on $\mathbb R^7$ discussed above. Moreover, there is a way to canonically determine both a metric and an orientation in a \emph{highly non-linear way} from the $3$-form $\varphi$. Then one can define a cross product $\times$ by using the metric to ``raise an index'' on $\varphi$. In summary, a manifold $(M, \varphi)$ with $\mathrm{G}_2$-structure comes equipped with a metric, cross product, $3$-form, and orientation that satisfy
\begin{equation*}
\varphi(u,v,w) = \langle u \times v , w \rangle.
\end{equation*}
This is exactly analogous to the data of an \emph{almost Hermitian manifold}, which comes with a metric, an almost complex structure $J$, a $2$-form $\omega$, and an orientation that satisfy
\begin{equation*}
\omega(u,v) = \langle Ju , v \rangle.
\end{equation*}
Essentially, a manifold admits a $\mathrm{G}_2$-structure if we can identify each of its tangent spaces with the imaginary octonions $\mathrm{Im} (\mathbb O) \cong \mathbb R^7$ in a smoothly varying way, just as an almost Hermitian manifold is one in which we can identify each of its tangent spaces with $\mathbb C^m$ (together with its Euclidean inner product) in a smoothly varying way. Manifolds with $\mathrm{G}_2$-structure also admit distinguished classes of \emph{calibrated submanifolds} similar to the pseudo-holomorphic curves of almost Hermitian manifolds. See~\cite{HL} for more about calibrated submanifolds.

For a manifold to admit a $\mathrm{G}_2$-structure, necessary and sufficient conditions are that it be \emph{orientable} and \emph{spin}, equivalent to the vanishing of its first two Stiefel-Whitney classes. Hence there are \emph{lots} of such $7$-manifolds, just as there are lots of almost Hermitian manifolds. But the story does not end there.

Let $(M, \varphi)$ be a manifold with $\mathrm{G}_2$-structure. Since it determines a Riemannian metric $g_{\varphi}$, there is an induced Levi-Civita covariant derivative $\nabla$, and one can ask if $\nabla \varphi = 0$. If this is the case, $(M, \varphi)$ is called a $\mathrm{G}_2$-manifold, and one can show that the Riemannian holonomy of $g_{\varphi}$ is contained in the group $\mathrm{G}_2 \subset \mathrm{SO}(7)$. Finding such ``parallel'' $\mathrm{G}_2$-structures is \emph{very hard}, because one must solve a fully non-linear partial differential equation for the unknown $3$-form $\varphi$. The $\mathrm{G}_2$-manifolds are in some ways analogous to \emph{K\"ahler} manifolds, which are exactly those almost Hermitian manifolds that satisfy $\nabla \omega = 0$. K\"ahler manifolds are much easier to find, partly because the metric $g$ and the almost complex structure $J$ on an almost Hermitian manifold are essentially independent (they just have to satisfy a mild condition of compatibility) whereas in the $\mathrm{G}_2$ case, the metric and the cross product are both determined non-linearly from $\varphi$. However, the analogy is not perfect, because one can show that when $\nabla \varphi = 0$, the Ricci curvature of $g_{\varphi}$ necessarily vanishes. So $\mathrm{G}_2$-manifolds are always \emph{Ricci-flat}! (This is one reason that physicists are interested in such manifolds---they play a role as ``compactifications'' in $11$-dimensional $M$-theory analogous to the role of Calabi-Yau $3$-folds in $10$-dimensional string theory. See~\cite{G} for a survey of the role of $\mathrm{G}_2$-manifolds in physics.) Thus in some sense $\mathrm{G}_2$-manifolds are more like \emph{Ricci-flat K\"ahler} manifolds, which are the \emph{Calabi-Yau} manifolds.

In fact, if we allow the holonomy to be a proper subgroup of $\mathrm{G}_2$, there are many examples of $\mathrm{G}_2$-manifolds. For example, the flat torus $T^7$, or the product manifolds $T^3 \times \mathrm{CY}2$ or $S^1 \times \mathrm{CY}3$, where $\mathrm{CY}n$ is a Calabi-Yau $n$-fold, have holonomy groups properly contained in $\mathrm{G}_2$. These are in some sense ``trivial'' examples because they reduce to lower-dimensional constructions. Manifolds with \emph{full holonomy} $\mathrm{G}_2$ are also called \emph{irreducible} $\mathrm{G}_2$-manifolds and it is precisely these manifolds that Bryant, Bryant--Salamon, Joyce, and Kovalev constructed.

We are still lacking a ``Calabi-Yau type'' theorem which would give necessary and sufficient conditions for a compact $7$-manifold that admits $\mathrm{G}_2$-structures to admit a $\mathrm{G}_2$-structure that is parallel ($\nabla \varphi = 0$.) Indeed, we don't even know what the conjecture should be. Some \emph{topological obstructions} are known, but we are far from knowing sufficient conditions. In fact, rather than comparing it to the Calabi conjecture, we should instead compare it to a different problem that it resembles more closely, namely, the following. Suppose $M^{2n}$ is a compact, smooth, $2n$-dimensional manifold that admits almost complex structures. What are necessary \emph{and sufficient} conditions for $M$ to admit K\"ahler metrics? We certainly know many necessary topological conditions, but we are nowhere near knowing sufficient conditions.

What makes the Calabi conjecture tractable (although certainly difficult) is the fact that we already start with a K\"ahler manifold (holonomy $\mathrm{U}(m)$ metric) and try to reduce the holonomy by only one dimension, to $\mathrm{SU}(m)$. Then the $\partial \bar \partial$-lemma in K\"ahler geometry allows us to reduce the Calabi conjecture to an (albeit fully non-linear) elliptic PDE for a \emph{single scalar function}. Any analogous ``conjecture'' in either the K\"ahler or the $\mathrm{G}_2$ cases would have to involve a \emph{system} of PDEs, which are much more difficult to deal with.

\subsection*{Acknowledgements}

An earlier version of this article appeared as an answer to a question on MathOverflow, and can be found here: \url{mathoverflow.net/questions/49357/g-2-and-geometry}. The author thanks Pete L. Clark for suggesting that it be suitable for a ``What is...'' article for the \emph{Notices of the AMS}, and also thanks the referee and editor for their useful suggestions.

\end{document}